\documentclass[a4paper,12pt]{elsarticle}
\usepackage{graphicx}
\usepackage{amsmath,mathtools}
\usepackage{parskip}
\usepackage{paralist}
\usepackage{amssymb}
\usepackage{pdfpages}
\usepackage{subcaption}
\usepackage[labelfont=bf]{caption}
\usepackage{sidecap}
\usepackage[titletoc,toc,title]{appendix}
\usepackage[margin=2cm]{geometry}
\usepackage{float}
\usepackage{epstopdf}
\usepackage{tikz}
\usepackage{tkz-fct}
\usetikzlibrary{arrows,matrix,positioning,decorations.pathreplacing,shapes,calc,decorations.markings,math,intersections}
\usepackage{tkz-euclide}
\usepackage{tikz-3dplot}
\usepackage{pgfplots}
\pgfplotsset{compat=1.14}
\usepgfplotslibrary{fillbetween}
\usepackage{siunitx}
\usepackage{lineno}
\usepackage{natbib}
\usepackage[hidelinks]{hyperref}

\AtBeginDocument{}

\biboptions{comma,round,semicolon,authoryear}
\setcitestyle{notesep={; },authoryear,round,semicolon,open={(},close={)}}

\newcommand{\wordeqn}[3]{{\small \begin{bmatrix}\text{#1}\\ \text{#2} \\ \text{#3}\end{bmatrix} }}
\newcommand{\swordeqn}[2]{{\small \begin{bmatrix}\text{#1}\\ \text{#2} \end{bmatrix} }}

\newcommand{\SR}{ASR}
\newcommand{\srlong}{adult sex ratio}
\newcommand{\Srlong}{Adult sex ratio}

\tikzset{
    >=stealth',
    boxed/.style={
           rectangle,
           rounded corners,
           draw=black, very thick,
           text width=6.5em,
           minimum height=2em,
           text centered},
    arr/.style={
           ->,
           thick,
           shorten <=2pt,
           shorten >=2pt,}
}

\newenvironment{nalign}{
    \begin{equation}
    \begin{aligned}
}{
    \end{aligned}
    \end{equation}
    \ignorespacesafterend
}

\journal{Theoretical Population Biology}
\title{\Srlong\ as an index for male strategy in primates}
\author[1]{Danya Rose}
\author[2]{Kristen Hawkes}
\author[1]{Peter S. Kim}
\address[1]{School of Mathematics and Statistics, University of Sydney, Sydney, NSW 2006, Australia}
\address[2]{Department of Anthropology, University of Utah, Salt Lake City, UT 84112, USA}

\begin{document}
\begin{frontmatter}

\begin{abstract}
The \srlong\ (\SR) is defined as the number of fertile males divided by the number of fertile females in a population. We build an ODE model with minimal age structure, in which males compete for paternities using either a multiple-mating or searching-then-guarding strategy, to investigate the value of \SR\ as an index for predicting which strategy males will adopt, with a focus in our investigation on the differences of strategy choice between chimpanzees (\emph{Pan troglodytes}) and human hunter-gatherers (\emph{Homo sapiens}). Parameters in the model characterise aspects of life history and behaviour, and determine both dominant strategy and the \SR\ when the population is at or near equilibrium. Sensitivity analysis on the model parameters informs us that \SR\ is strongly influenced by parameters characterising life history, while dominant strategy is affected most strongly by the effectiveness of guarding (average length of time a guarded pair persists, and resistance to paternity theft) and moderately by some life history traits. For fixed effectiveness of guarding and other parameters, dominant strategy tends to change from multiple mating to guarding along a curve that aligns well with a contour of constant \SR, under variation of parameters such as longevity and age female fertility ends. This confirms the hypothesis that \SR\ may be a useful index for predicting the optimal male mating strategy, provided we have some limited information about ecology and behaviour.
\end{abstract}
\end{frontmatter}

\section{Introduction}
\label{sec:intro}
The closest living genus to our human genus, \emph{Homo}, is genus \emph{Pan}; and we share many physiological, developmental and behavioural traits with them. Significant differences exist, however, particularly regarding life history and the social structure around mating arrangements. While both humans and chimpanzees engage in a variety of strategies, such as multiple mating, possessive short-term, or longer-lasting exclusive relationships (see, for example, \citep{tutin1979}), each species tends to engage in one class of strategies with greater frequency than the others. Some recent studies examine the evolution of monogamy from a mathematical perspective; for example, \citep{looChanHawkesKim2017,looHawkesKim2017,schachtBell2016,schachtKramerSzekelyKappeler2017} all discuss the role of mating sex ratio and partner availability in the evolution of monogamy or other mating-related behaviour. In our study, we broadly categorise the reproductive strategies as either multiple mating or mate guarding, in order to build a relatively simple model that captures sufficient dynamics to explore the problem of predicting strategy by observing demography. Specifically, chimpanzees (with relatively more females per male) typically engage in multiple mating more frequently than guarding, and hunter-gatherers (with relatively scarce fertile females) tend to engage in guarding more frequently than multiple mating.

Life tables indicate that roughly half of infants in wild chimpanzee and hunter-gatherer groups will die before reaching maturity (sexual maturity). On reaching maturity, though, chimpanzees and humans take divergent paths: a chimpanzee just at maturity (first birth around age~$14$) may expect to live another fifteen to twenty years \citep{hillEtAl2001}, but a human hunter-gatherer at the same point (first birth around age~$19$) may expect to live healthily for another forty years or more \citep{blurtonJones2016,hillHurtado1995,howell1979}. The age of last birth occurs in both humans and chimps at around~$45$ years, and almost no wild chimpanzee can expect to survive until this age; those who do, though robust enough to have reached that age, are by that time relatively frail. In stark contrast, many hunter-gatherer females live in good health beyond menopause, only becoming frail into their seventies; those who survive so long have spent half their adult lives post-fertile.

Though pair bonding is not unknown amongst the primates, humans are unique among great apes in that we tend to form relatively long-lasting pair bonds, as opposed to multiple mating, wherein no long-term attachments are made. A favoured hypothesis to explain our different behaviour is that it arises from paternal investment \citep{washburnlancaster1968,lancasterLancaster1983}---the hypothesis that care provided by males (ostensibly in the form of meat acquired through their hunting) to their offspring is the basis for the formation of cooperative pair bonds between men and women to help raise their young---but more recent work suggests that other mechanisms may be responsible for this difference in behaviour from our relatives (see, for example, \citep{gurvenHill2009,hawkesEtAl2010,looChanHawkesKim2017,looHawkesKim2017,lukasCluttonbrock2013}). 
The extension of the human life span (and of the male fertility period) without lengthening female fertility to later ages directly changes the ratio of fertile males to fertile females (called the \srlong: \SR). As the \SR\ increases, the number of fertile females available per male decreases, which we argue changes the incentives for adopting either strategy.

Female choice arises as an important question in mating dynamics, and incorporating that dynamic in game-theoretical models of mating strategy selection can significantly change the outcome of such a model, as discussed in, for example, \citep{parkerBirkhead2013,cluttonbrockParker1995,opieEtAl2013}. As the focus of our study is the chimpanzee-human relationship, we make the explicit decision to disregard female choice, for the following reasons. Chimpanzee males tend to coerce---often violently---fertile females to have sex, potentially resulting in multiple partners for females \citep{mullerEtAl2007}. Furthermore, in many primates the risk of infanticide by males is very high \citep{vanschaikCarelJanson2000} (including among chimpanzees: see \citep{arcadiWrangham1999} and references therein), leading to varied strategies, including female promiscuity to confuse paternity \citep{mullerEtAl2007}; consequently, females chimpanzees are not especially choosy. In humans, even, choosiness among females is not necessarily optimal for the survival of their offspring, as discussed in \citep{hrdy2003}, for instance, leading to polyandrous mating.

\section{Methods}
\label{sec:methods}
In this paper, we construct a simple two-strategy ODE model, in which males either guard mates (once acquired) or multiply mate (that is, possibly acquiring many mates at the same time), and competition between strategies occurs only through the acquisition of paternities. This approach is similar to the dynamic analysis assessed in \citep{mylius1999}, which explicitly models pair formation. Our approach, like the work of \citep{kvrivanCressman2017,kvrivanGalanthayCressman2018}, explicitly acknowledges finite (and differing) interaction times affecting the payoff (in the form of paternities) for the different strategies, without the use of traditional game-theoretical payoff matrices.

Offspring are assumed to inherit the strategies of their fathers. Guarding males do not contribute to the survival of their offspring, and females do not express preference for either type of male. Our model has a number of parameters, which correspond to aspects of life history, ecology, and behaviour. We will vary the parameters and record the resulting dominant male strategy and the resulting \SR, with the intent to determine whether, or under what circumstances, the \SR\ can serve as an index for determining what strategy males are most likely to employ.

Inheritance of strategy by paternal descent is assumed as a method of simplification. It is of course very likely that multiple factors contribute to the strategy that a male may choose, and that the strategy an individual carries may be ``mixed'' (in the sense that he commits some fraction of his efforts to one strategy at the expense of full efficiency in the other), and in reality an individual may conceivably change his strategy investment over the course of his own life. Our simplification ensures that the representation of each strategy across a population is in proportion from one generation to the next.

\subsection{Life cycle}
\label{sec:lifeCycle}
\begin{figure}\centering
\begin{tikzpicture}[node distance=2.5cm, auto,rectangle,rounded corners,line width=0.2mm]
	\node[draw]			(F)		at( 0, 0)		{$F$};
	\node[draw]			(G)		at(-3, 0)		{$G$};
	\node[draw]			(M)		at( 3, 0)		{$M$};
	\node[draw]			(FG)	at(-1,-4)		{$F^G$};
	\node[draw]			(FM)	at( 1,-4)		{$F^M$};
	\node[draw]			(CG)	at(-5,-2)		{$C^G$};
	\node[draw]			(CM)	at( 5,-2)		{$C^M$};
	\node[align=center]	(FMM)	at( 1.3,-1.8)	{Mating\\occurs};
	\node[align=center]	(FGG)	at(-1.3,-1.8)	{Pair\\forms};
                        
	\draw[->]			(F)		to	(FGG);
	\draw[->]			(G)		to	(FGG);
	\draw[->]			(FGG)	to	(FG);
                        
	\draw[->]			(F)		to	(FMM);
	\draw[->]			(FMM)	to	(FM);
	\draw[<->]			(M)		to	(FMM);
                        
	\draw[->,dashed]	(FG)	to	(CG);
	\draw[->,dashed]	(FM)	to	(CM);
                        
	\draw[->]			(CM)	to	(M);
	\draw[->,dashed]	(FG)	to	(CM);
                        
	\draw[->]			(CG)	to	(F);
	\draw[->]			(CG)	to	(G);
                        
	\draw[->]			(CM)	to	(F);
                        
	\draw[->]			(FG)	to	(G);
	\draw[->]			(FG)	to	(F);
                        
	\draw[->]			(FM)	to	(F);
\end{tikzpicture}
\caption{Illustration of our model, excluding removal due to death (adults at rate $\mu$ per head, young at rate $\delta$ per head) or retirement (adults only at rates $\tau$ for females and $\lambda$ for males). Units of $F^G$ are pairs of individuals; if either partner dies or retires the remaining partner returns to the original group. Multiple-mating males will mate with receptive females and return immediately to search for another receptive female. Some offspring of guarded pairs have multiple-mating fathers due to paternity theft.}
\label{fig:model}
\end{figure}
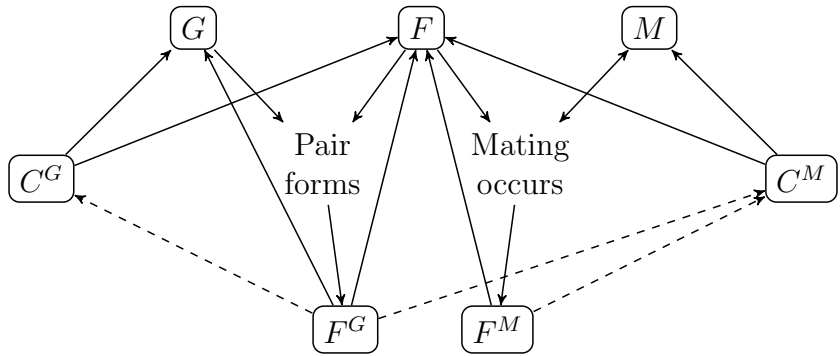

Consider a population $P$ comprising
\begin{enumerate}
\item fertile searching males $G$, who guard their mates once paired,
\item fertile multiple-mating males $M$, who return to the mating pool as soon as they have mated,
\item fertile receptive females $F$, who may be ``recruited'' by a $G$ or $M$ male and may or may not have been previously recruited,
\item fertile unreceptive females $F^M$, who have been recruited by a multiple-mating male and are thus occupied short term by that male, or pregnant or lactating as a consequence,
\item guarded pairs $F^G$, consisting of one fertile guarding male whose energy goes into (possibly imperfectly) guarding his mate from multiple maters and one fertile female who may bear offspring to multiple maters who steal paternities from her partner,
\item offspring of guarding males $C^G$ (called ``guarded offspring''), and
\item offspring of multiple-mating males $C^M$ (called ``unguarded offspring'', though may be born to guarded females),
\end{enumerate}
summarised in Table~\ref{tab:vars}. The total population is $P = F + G + M + 2F^G + F^M + C^G + C^M$, where the number of guarded pairs $F^G$ is counted twice because each guarded pair contains two individuals. Parameters of the model are summarised in Table~\ref{tab:params}, and its dynamics illustrated in Figure~\ref{fig:model}, and all are described below.
\begin{table}[b]
\centering
\caption{List and descriptions of all variables. These variables represent number, rather than proportion of population.}
\label{tab:vars}
\begin{tabular}{cl}\hline
Variable	& Description					\\\hline
$P$			& Total population				\\
$F$			& Receptive females				\\
$G$			& Searching (guarding) males	\\
$M$			& Multiple-mating males			\\
$F^G$		& Guarded pairs					\\
$F^M$		& Unreceptive females 			\\
$C^G$		& Guarded offspring				\\
$C^M$		& Unguarded offspring			\\\hline
\end{tabular}
\end{table}

Juveniles are assumed to mature at rate $\gamma$ per individual, and die at rate $\delta$ per individual (assumed independent of population density). Half the juveniles who survive to maturity we assume\footnote{Although there is variability of birth sex ratios among primates in general---in some instances quite far from $0.5$---it remains that in the majority of cases the ratio is close to $0.5$. A table of birth sex ratios for $102$ species is given in \citep[][supplemental material]{silkBrown2008}, for instance, of which $42$ also include wild data.} become receptive females, while the other half become either multiple maters or searchers (respectively the same as their fathers). Observe that paternal investment is not modelled here. Unguarded offspring do not receive less care than guarded offspring, and no distinction is made among the survival of juveniles depending on their paternity.

Fertile females include receptive and unreceptive females, and female members of guarded pairs. They die at base rate $\mu$ (assumed to be independent of population density) and advance to menopause at rate $\tau$. Post-fertile females ``retire'' from the population as they reach menopause, as they no longer play a role in the reproductive dynamics in this model. We do not explicitly model grandmother effects, such as extra calorific benefits due to the presence of post-fertile females who can provision the offspring of other females (see, for instance, \cite{coxworthKimMcQueenHawkes2015,looHawkesKim2017} for connections between grandmothering and sex-ratio, and \cite{kimCoxworthHawkes2012,kimMcQueenCoxworthHawkes2014,kimMcQueenHawkes2019} for modelling of the relationship between grandmothering and human life history). Fertile males include searching males, multiple-mating males and the male partners in guarded pairs. Fertile males have a base death rate of $k\mu$, and retire due to geriatric infirmity at rate $\lambda$, whereupon they too leave the population.

Receptive females are ``recruited'' by either searching or multiple-mating males at rate $r$ per possible pair. If recruited by a searching male, they together form a guarded pair, which will on average break spontaneously with rate $\beta$ per pair, with both partners returning to their respective original pools. Guarded pairs may also be broken by the death or retirement of one partner, in which case the remaining partner returns to his or her own original pool. Females recruited by multiple-mating males become unreceptive for a period of time, returning to the receptive pool with average rate $\sigma$. Females in guarded pairs and unreceptive females produce offspring (respectively guarded and unguarded offspring) at rate $\rho$ per bonded/unreceptive female, which includes time spent pregnant or lactating. We can think of $\frac{1}{\rho}$ as the inter-birth interval. Note that in our model females in the $F^G$ and $F^M$ pools are \emph{producing} offspring (as opposed to necessarily having current dependents), and even if they retire or die, we assume that any juveniles already born are cared for until maturity or death. That is, our model assumes that all care required by a juvenile is provided, no matter how much or how little that is.

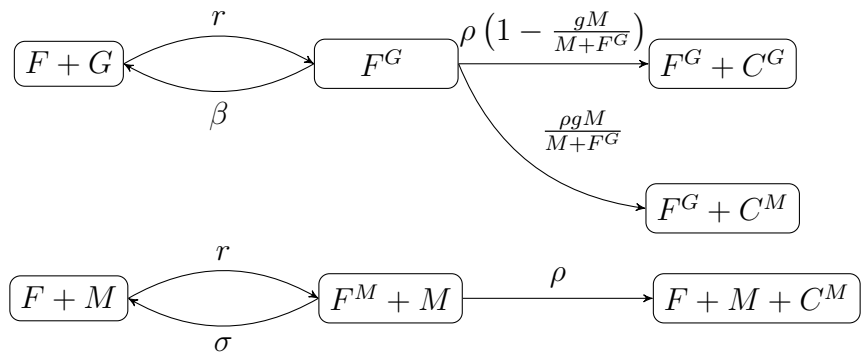
\begin{figure}\centering
\begin{tikzpicture}[node distance=2.5cm, auto,rectangle,rounded corners]
	\node[draw]								(FPG)		{$F+G$};
    \node[draw,right=of FPG]				(FG)		{\phantom{$+$} $F^G$ \phantom{$M$}};
    \node[draw,right=of FG]					(FGCG)		{$F^G+C^G$};
    \node[draw,below= 1.25cm of FGCG]		(FGCM)		{$F^G+C^M$};
    \draw[->,bend left]		(FPG.east)	to node[above,pos=0.5]			{$r$}			(FG.west);
    \draw[->,bend left]		(FG.west)	to node[below,pos=0.5]			{$\beta$}		(FPG.east);
    \draw[->]				(FG.east)	to node[above,pos=0.5]			{$\rho\left(1-\frac{gM}{M+F^G}\right)$}	(FGCG.west);
    \draw[->,bend right]	(FG.east)	to node[above right,pos=0.5]	{$\frac{\rho gM}{M+F^G}$}		(FGCM.west);
    \node[draw,below=of FPG]				(FPM)		{$F+M$};
    \node[draw,right=of FPM]				(FMPM)		{$F^M+M$};
    \node[draw,right=of FMPM]				(FMPMPCM)	{$F+M+C^M$};
    \draw[->,bend left]		(FPM.east)	to node[above,pos=0.5]			{$r$}			(FMPM.west);
    \draw[->,bend left]		(FMPM.west)	to node[below,pos=0.5]			{$\sigma$}		(FPM.east);
    \draw[->]				(FMPM.east)	to node[above,pos=0.5]			{$\rho$}		(FMPMPCM.west);
\end{tikzpicture}
\caption{Illustration of paternity theft dynamics. Receptive females $F$ may pair with either guarding $G$ or multiple-mating $M$ males; females who have mated with a male of either type produce offspring at rate $\rho$, but a proportion of pair-bonded females produce offspring fathered by multiple maters.}
\label{fig:paternityTheft}
\end{figure}

Once a searching guarder has recruited a receptive female, he stops searching and invests his energy in guarding, so as to be assured of the paternity of the offspring his mate produces, and stops seeking other mates. Multiple-mating males, however, return immediately to the same pool to search for other mates. The kinetics of this mating model are as follows:
\begin{equation*}
\begin{array}{clccl}
F + G & \xrightleftharpoons[\beta] {\hspace{0.5cm}r\hspace{0.5cm}} & F^G     &\xrightarrow{\hspace{0.5cm}\rho\hspace{0.5cm}} & F^G + C^G, \\[1em]
F + M & \xrightleftharpoons[\sigma]{\hspace{0.5cm}r\hspace{0.5cm}} & F^M + M &\xrightarrow{\hspace{0.5cm}\rho\hspace{0.5cm}} & F + M + C^M.
\end{array}
\end{equation*}
However, multiple maters may attempt to steal paternities from males in guarded pairs, and succeed at rate $g$ with respect to the density of multiple maters, illustrated in Figure~\ref{fig:paternityTheft}. Paternity theft is taken to be in proportion to $M/(M+F^G)$, the density of multiple-mating males amongst all males who are considered to be in line for paternity, thus representing a frequency-dependent interaction rate. Searching $G$ males are excluded from the density expression because although they are searching for mates, if they obtain mates from the pool of receptive females $F$ they enter the pool $F^G$ of guarded pairs, but if they \emph{were} to steal a paternity from males already in guarded pairs, according to our strategy inheritance model the resulting offspring will still be guarded $C^G$, the male fraction of whom will mature into the pool of searching guarding males $G$; hence paternity theft by guarders is not considered.

Together, $\beta$ and $g$ constitute an ``effectiveness of guarding'': $\beta$ relates to the ability of guarding males to retain their mates for extended periods of time, and $g$ relates to their ability to ensure that the offspring born to their female partners are indeed their own.

Finally, we include a population-density-dependent death rate $\nu P$, so that population growth is constrained. Removal due to the population density dependence does not have to be considered as death, but could also represent migration away from the population of interest. At this point, we can write the equations. A dot denotes a time derivative.
\begin{nalign}\label{eqn:DEsystem}
\dot F   &= \frac{1}{2}\gamma \left(C^G + C^M\right) + \left(\beta + \lambda + \mu \left(k+\nu P\right)\right) F^G + \sigma F^M - \left(r\left(G + M\right) + \tau + \mu \left(1+\nu P\right)\right) F, \\
\dot G   &= \frac{1}{2}\gamma C^G + \left(\beta + \tau + \mu \left(1+\nu P\right)\right) F^G - \left(rF + \lambda + \mu \left(k+\nu P\right)\right) G, \\
\dot M   &= \frac{1}{2}\gamma C^M - \left(\lambda + \mu \left(k+\nu P\right)\right) M, \\
\dot F^G &= rFG - \left(\beta + \tau + \lambda + \mu \left(1+k+2\nu P\right)\right) F^G, \\
\dot F^M &= rFM - \left(\sigma + \tau + \mu \left(1+\nu P\right)\right) F^M, \\
\dot C^G &= \rho\left(1 - g\frac{M}{M+F^G}\right) F^G  - \left(\gamma + \delta \left(1+\nu P\right)\right) C^G, \\ 
\dot C^M &= \rho\left(F^M + g\frac{MF^G}{M+F^G}\right) - \left(\gamma + \delta \left(1+\nu P\right)\right) C^M.
\end{nalign}
In words, this is
\begin{align*}
\dot F   &= \swordeqn{rate female}{juveniles mature} + \wordeqn{rate guarded pairs}{break spontaneously or}{by death or retirement}+\swordeqn{rate unreceptive}{females return}\\&\qquad - \left(\wordeqn{rate receptive}{females recruited}{by males} + \wordeqn{rate receptive}{ females retire}{or die}\right), \\
\dot G   &= \swordeqn{rate guarded male}{offspring mature}+\wordeqn{rate guarded pairs}{break spontaneously or}{by death or retirement}-\left(\swordeqn{rate guarded}{pairs form}+\wordeqn{rate guarding}{males retire}{or die}\right), \\
\dot M   &= \swordeqn{rate unguarded male}{offspring mature}-\wordeqn{rate multiple-}{maters retire}{or die}, \\
\dot F^G &= \swordeqn{rate guarded}{pairs form}-\swordeqn{rate guarded}{pairs break}, \\
\dot F^M &= \wordeqn{rate receptive}{females recruited}{by multiple maters}-\wordeqn{rate unreceptive}{females return,}{retire or die}, \\
\dot C^G &= \swordeqn{rate guarded}{offspring are born}-\swordeqn{rate guarded}{offspring mature or die}, \\
\dot C^M &= \wordeqn{rate unguarded}{offspring are born,}{including paternity theft}-\swordeqn{rate unguarded}{offspring mature or die}.
\end{align*}
Offspring of multiple maters born to guarded females are automatically classified as unguarded (due to the terms $\pm gMF^G/(M+F^G)$ in the last two lines of Equation~\eqref{eqn:DEsystem} above, indicating paternity theft), because strategies are assumed for simplicity to be inherited patrilineally. This means a reduction in the rate at which guarded offspring are born and a corresponding increase in the rate at which unguarded offspring are born when both $g$ and $M$ are not zero.

\subsection{Life history model}
\label{sec:model}
Our principal interest is in understanding the relative success of alternative strategies depending on fertility and mortality parameters, and guarding effectiveness. 
Mortality rates at all ages determine the average longevities, but we are interested in what strategy (or balance of strategies) might be chosen by populations with given longevities. Dividing the population into fertile adults and juveniles imposes a minimal age structure, which we can conveniently use to reframe the parameters. We wish to exchange $\mu$, $\delta$ and $\gamma$ for parameters that characterise certain aspects of life history and environment.

Our model takes advantage of the elegance and simplicity of ODEs to achieve its goals. As such, it is not the same as an age-structured PDE model, although it contains a basic age structure which we shall use. We have two compartments representing juveniles and five compartments representing adults. Individuals begin in a juvenile compartment and move to an adult compartment at some time in the interval $(0,\infty)$ to take part in the adult dynamics that are of interest. As a result, some individuals can transition into maturity at arbitrarily low age, and some individuals can remain juvenile for an arbitrarily long time; in a certain sense, an ODE system such as this can only speak in terms of probabilities and of average times. Later work may introduce explicit age structure in an attempt to build in more realism.

\begin{figure}
\begin{center}
\begin{tikzpicture}[node distance=2.5cm, minimum width=1.5cm, minimum height=0.5cm, auto,rectangle,sharp corners]
	\node[draw]								(C)		{$C$};
    \node[draw,right=of C]					(A)		{$A$};
	\node[below=1.5cm of C]					(CD)	{};
	\node[below=1.5cm of A]					(AD)	{};
	\draw[->]				(C)	to node[above] {$\gamma$}	(A);
	\draw[->]				(C)	to node[right] {$\delta$}	(CD);
	\draw[->]				(A)	to node[right] {$\mu$}		(AD);
\end{tikzpicture}
\end{center}
\caption{Block diagram of the age structure model governed by Equations~\eqref{eqn:CADEs}, with juveniles denoted by $C$ and adults by $A$. Juvenile death rate, adult death rate and maturation rate are represented by $\delta$, $\mu$ and $\gamma$ respectively.}
\label{fig:ageStructureBlock}
\end{figure}
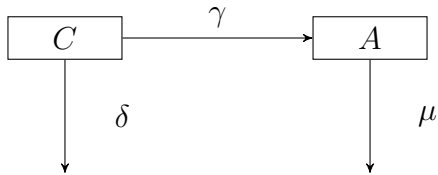

If we consider only females and combine all female juveniles and all female adults into respective groups $C$ and $A$, and consider only the dynamics of maturation (due to $\gamma$) and the linear death rates $\delta$ and $\mu$ independent of population density (otherwise due to $\nu$), we can extract dynamics that represent, in some sense, the female survivorship and the demographic distribution (between juveniles and adults) of the female population of age $t$. This maturation and death dynamic is illustrated in Figure~\ref{fig:ageStructureBlock}, where the transition from $C$ to $A$ is governed by a intensity $\gamma$, rather than a strict time (as in a model with more detailed ageing dynamics). These dynamics will be governed by differential equations \begin{nalign}\label{eqn:CADEs}\frac{dC}{dt} &= -(\delta+\gamma)C, \\ \frac{dA}{dt} &= \gamma C - \mu A,\end{nalign} and we will take initial conditions $C(0) = 1$ and $A(0) = 0$. The solution of equations \eqref{eqn:CADEs} is easily obtained by standard techniques (e.g. elimination of one variable and the use of an integrating factor), yielding two exponentially decreasing functions, for positive $\delta$, $\gamma$ and $\mu$.

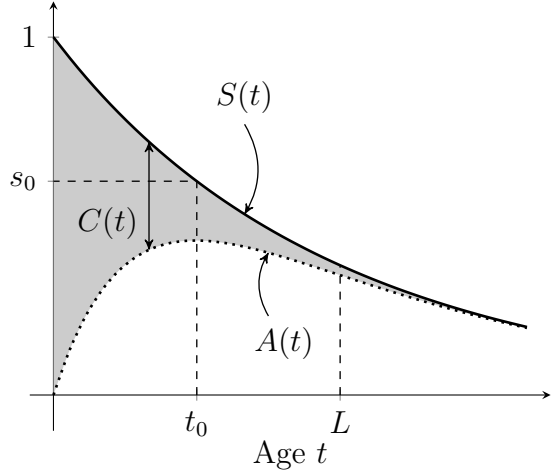
\begin{figure}
\begin{center}
\newcommand*{\deltaval}{0.017571954578652}%
\newcommand*{\gammaval}{0.041805946202065}%
\newcommand*{\muval}{0.016313407572677}%
\def\Lfun(#1,#2,#3){(#1+#3)/(#3*(#1+#2))}%
\def\tzerofun(#1,#2,#3){\Lfun(#1,#2,#3)/2}%
\pgfmathsetmacro{\tzeroval}{\tzerofun(\gammaval,\deltaval,\muval)}%
\pgfmathsetmacro{\ezeroval}{\Lfun(\gammaval,\deltaval,\muval)}%
\def\Afun(#1){(\gammaval*exp(-\muval*#1))/(\deltaval + \gammaval - \muval) - (\gammaval*exp(-#1*(\deltaval + \gammaval)))/(\deltaval + \gammaval - \muval)}%
\def\Sfun(#1){(\gammaval*exp(-\muval*#1))/(\deltaval + \gammaval - \muval) - (\gammaval*exp(-#1*(\deltaval + \gammaval)))/(\deltaval + \gammaval - \muval) + exp(-#1*(\deltaval + \gammaval))}%
\pgfmathsetmacro{\szeroval}{\Sfun(\tzeroval)}%
\begin{tikzpicture}
\begin{axis}[domain=0:100,xmax=105,ymax=1.1,xmin=-5,ymin=-0.1,samples=30, axis x line=middle, axis y line=center, xtick={0,\tzeroval,\ezeroval}, xticklabels={$0$,$t_0$,$L$}, ytick={0,\szeroval,1}, yticklabels={$0$,$s_0$,$1$}, xlabel= Age $t$, x label style={at={(axis description cs:0.5,0)},anchor=north},]
	\addplot[name path=line,smooth,solid,black,line width=1] {\Sfun(x)};
	\addplot[name path=dotted,smooth,dotted,black,line width=1] {\Afun(x)};
	\addplot[gray!40] fill between[of=dotted and line,];
	\draw[dashed,black,line width=0.5]	(\tzeroval,0) -- (\tzeroval,\szeroval) -- (0,\szeroval);
	\draw[dashed,black,line width=0.5]	(\ezeroval,0) -- (\ezeroval,{\Sfun(\ezeroval)});
	\draw[->,black,line width=0.5] ({(\tzeroval+\ezeroval)/2},{\Afun((\tzeroval+\ezeroval)/2)/3}) node[below] {$\quad A(t)$} to[bend left] ({((\tzeroval+\ezeroval)/2)},{\Afun({((\tzeroval+\ezeroval)/2)})});
	\draw[->,black,line width=0.5] ({(\tzeroval+\tzeroval+\ezeroval)/3.0},{1.5*\Sfun({((\tzeroval+\tzeroval+\ezeroval)/3.0)})}) node[above] {$S(t)$} to[bend left] ({((\tzeroval+\tzeroval+\ezeroval)/3.0)},{\Sfun({((\tzeroval+\tzeroval+\ezeroval)/3.0)})});
	\draw[<->,black,line width=0.5] ({\tzeroval*2/3},{\Afun(\tzeroval*2/3)}) node[above left] {$C(t)$} -- ({\tzeroval*2/3},{\Sfun(\tzeroval*2/3)});
\end{axis}
\end{tikzpicture}
\end{center}
\caption{A typical example of $S(t)$ (solid curve) governed by the dynamics of equations~\eqref{eqn:CADEs}. The dotted curve shows $A(t)$, indicating the demographic breakdown of individuals of age $t$, split between juveniles and adults. The grey shaded region indicates juveniles: initially all individuals are juveniles, who decrease in proportion of surviving population exponentially as age $t\to\infty$ according to $\exp(-(\delta+\gamma)t)$ (through either transition into maturity or death). At age $t_0=L/2$ we have proportion $s_0 = S(t_0)$ of initial population surviving.}
\label{fig:survivalCurve}
\end{figure}

Because we take the initial conditions of Equations~\eqref{eqn:CADEs} such that $C(0)+A(0) = 1$, the sum $S(t) := C(t)+A(t)$ represents the probability of an individual surviving to age $t$, depicted in Figure~\ref{fig:survivalCurve} as the solid curve. Note that $S(t)$ is a strictly decreasing function for $t\ge0$, whose limit is $0$ as $t\to\infty$. The ratio $\frac{A(t)}{S(t)}$ is the probability that a surviving individual has matured by age $t$; it is a simple exercise to check that this expression has a limiting value of $1$ as $t\to\infty$ when $\delta+\gamma>\mu$,\footnote{The case where $\delta+\gamma<\mu$ is not of interest to our particular study, as adult mortality is generally lower than juvenile mortality, but in that case the limiting ratio of adults over adults and juveniles is $\frac{\gamma}{\mu-\delta}$.} and by construction the mean transition time to maturity is $\frac{1}{\gamma}$. Also depicted in Figure~\ref{fig:survivalCurve} is $A(t)$ (dotted curve) and a representation of $C(t)$ as the vertical distance between the solid and dotted curves.

The expected time an individual lives in the population is $$L = \int_0^\infty S(t)\;dt = \frac{\gamma + \mu}{\mu\left(\delta + \gamma\right)}.$$ We impose the constraints $\gamma = 2/L = 1/t_0$, where $t_0$ is the mean time to reproductive maturity, and $S(t_0) = s_0$ (the proportion surviving at the mean transition time to maturity). The purpose of these constraints is to introduce parameters characterising certain aspects of the population: $s_0$ characterises something of the ``toughness'' of being a juvenile; and $L$ characterises the average life span of an individual.

Our constraints and the definition of $L$ give us a system of equations \begin{nalign}\label{eqn:paramsystem}s_0 &= S(t_0), \\ L &= \frac{\gamma + \mu}{\mu\left(\delta + \gamma\right)},\end{nalign} from which we can numerically obtain values of $\delta$ and $\mu$, given an expected life span $L$ and survival rate into maturity $s_0$.

As adult males are subject to base mortality rate $k\mu$, they suffer life expectancy at birth $\tilde L = \frac{\gamma + k\mu}{k\mu\left(\delta + \gamma\right)}$. Their relative life expectancy at birth is $\frac{\tilde L}{L} = \frac{\gamma + k\mu}{k\left(\gamma + \mu\right)}.$

Adult females have a fertile window of $t_1 - t_0$; we can thus set $\tau(L,t_1) = \frac{2}{2t_1 - L}$. Similarly, males have an fertile window of $t_2 - t_0$; and so we set $\lambda(L,t_2) = \frac{2}{2t_2 - L}$. For our purposes, we choose $t_1 = 45$ for the age of female fertility ends, and $t_2 = 75$ for the age of male retirement. While male chimpanzees cannot expect to live until such ages, such a high retirement age ensures that most low-longevity males will be removed by death, rather than retirement. In reality, it is likely that the age of male retirement is coupled somehow to longevity (as would be the age of female retirement/menopause in all primates but humans), but for the purposes of exploration in this model we keep these parameters decoupled.

The original parameters involved in life history ($\gamma$, $\delta$, $\mu$, $\tau$, $\lambda$) may now be replaced by expressions involving ($L$, $s_0$, $t_1$, $t_2$). Although $k$ affects the life span of males, we place it among the parameters we consider to represent behavioural aspects of the biology, ($r$, $g$, $\beta$, $\sigma$, $k$). The remaining two parameters, ($\rho$, $\nu$) we consider ecological/biological not related directly to life history or behaviour (as far as this study goes). Table~\ref{tab:params} lists and describes all parameters and typical values or ranges that we explore, and Table~\ref{tab:LHfuns} lists the mortality and other life history functions.

For Table~\ref{tab:params} we draw on \citep{blurtonJones1986,blurtonJones2016,emeryThompsonEtAl2007,gurvenKaplan2007,hawkesBlurtonJones2005,hillEtAl2001,nishidaEtAl2003} for life history data, which guides our choices for life history parameter values and ranges to explore. A noted shortcoming of our model is that realistic birth rates (denoted by $\rho$) tend to be slightly below replacement value; rather than compensate by increasing the mean window of female fertility (by decreasing either or both of the death and female retirement rates), we increase the birth rate to prevent extinction. Other parameters, such as population density-dependent death rate $\nu$, couple-forming rate $r$, paternity theft success rate $g$, pair-bond break-up rate $\beta$, and return rate $\sigma$ are guessed as being reasonable to explore. For instance, pair-bond stability is determined by $\beta$; human pair-bond length, for instance, is famously variable and as a topic consumes a large amount of glossy paper every year even for modern humans, yet established pairs sometimes last until both partners are senescent. Male-specific mortality in adulthood is typically higher, on average, than female (see, for example, \citep{courtenaySantow1989} for chimpanzees, or \citep{hillHurtadoWalker2007} for hunter-gatherers), though we choose to additionally explore a region where male-specific mortality is lower.

\begingroup\def\arraystretch{1.2}%
\begin{table}
\centering
\caption{List and descriptions of all parameters. Some parameters are chosen zero during initial analysis. Values and ranges for these parameters are chosen for relevance to understanding differences between hunter-gatherer and chimpanzee male mating strategies.}
\label{tab:params}\small
\begin{tabular}{cclp{4.5cm}}
\hline
Parameter			& Typical value/range						& Description													& Reference		\\
\hline
$L$					& $10$ to $50$								& Mean female longevity	& \citep{gurvenKaplan2007,hawkesBlurtonJones2005,hillEtAl2001}		\\
$s_0$				& $\tfrac{1}{3}$ to $\tfrac{2}{3}$			& Proportion of juveniles surviving				& \citep{gurvenKaplan2007,hillEtAl2001}		\\
$t_1$				& $30$ to $60$								& Age female fertility ends	& \citep{emeryThompsonEtAl2007,hawkesBlurtonJones2005,nishidaEtAl2003}		\\
$t_2$				& $60$ to $80$								& Age of male retirement	& \citep{gurvenKaplan2007,hawkesBlurtonJones2005,hillEtAl2001}	\\
$\rho$				& $\tfrac{1}{4}$ to $\tfrac{1}{2}$			& Rate juveniles are born	& \citep{blurtonJones1986,blurtonJones2016,nishidaEtAl2003}		\\
$\nu$				& $\tfrac{1}{1500}$ to $\tfrac{1}{500}$		& Crowding factor												& 		\\
$r$					& $\tfrac{1}{2}$ to $2$						& Couple-forming rate											& 		\\
$g$					& $0$ to $0.225$							& Paternity theft success rate									& 		\\
$\beta$				& $0$ to $\tfrac{1}{4}$						& Break-up rate for guarded pairs								& 		\\
$\sigma$			& $\tfrac{1}{2}$ to $2$						& Return rate for $F^M$ females									& 		\\
$k$					& $0.9$ to $1.1$							& Male death rate modifier										& \citep{hillEtAl2001}		\\
$\tau$				& see Table~\ref{tab:LHfuns}				& Rate of menopausal retirement									& 		\\
$\lambda$			& see Table~\ref{tab:LHfuns}				& Male retirement rate											& 		\\
\hline
\end{tabular}
\end{table}
\endgroup

Parameters can be grouped according to three broad themes: life history, ($L$, $s_0$, $t_1$, $t_2$); behaviour, ($r$, $g$, $\beta$, $\sigma$, $k$); and environment/biology, ($\rho$, $\nu$).

\begingroup\def\arraystretch{1.5}%
\begin{table}
\centering
\caption{Life history functions.}
\label{tab:LHfuns}
\begin{tabular}{ccl}
\hline
Parameter			& Functional form			& Description															\\
\hline
$t_0(L)$			& $L/2$						& Half mean life span													\\
$\gamma(L)$			& $2/L$						& Maturation rate of juveniles into adults	\\
$\delta(L,s_0)$		& Determined numerically from \eqref{eqn:paramsystem}	& Base juvenile death rate					\\
$\mu(L,s_0)$		& Determined numerically from \eqref{eqn:paramsystem}	& Base death rate of adults					\\
$\tau(L,t_1)$		& $\tfrac{1}{t_1-t^*(L)}$	& Transition rate of fertile females to post-fertile					\\
$\lambda(L,t_2)$	& $\tfrac{1}{t_2-t^*(L)}$	& Transition rate of fertile males to infertile\vspace{0.1em}			\\
\hline
\end{tabular}
\end{table}
\endgroup

\subsection{\Srlong\ (\SR)}
\label{subsec:SR}
The \srlong\ (\SR) is defined as the ratio of fertile males to fertile females. Our model assumes all adult males are fertile. Thus the \SR\ $a$ is given by $$a = \frac{\bar M}{\bar F},$$ where $\bar M = M+G+F^G$ is the total number of fertile males, and $\bar F = F+F^G+F^M$ is the total number of fertile females. The hypothesis of \citep{coxworthKimMcQueenHawkes2015} regarding humans, also the subject of \citep{looChanHawkesKim2017,looHawkesKim2017,schachtBell2016} is that the \SR\ determines male mating strategy. We aim to investigate whether the \SR\ is a sufficient index to determine the strategy (when there are two strategies in play) or if other parameters are necessary.

The model dynamics cause us to predict that increasing $t_1$ will decrease the \SR\ by causing females to remain fertile for longer, and thus increase relative to the number of fertile males. Similarly, increasing $t_2$ will increase the \SR\ by increasing the relative number of fertile males to fertile females. Increasing time that individuals are in the population $L$ (independently of $t_1$ and $t_2$) should increase the effect of disparity between $t_1$ and $t_2$, so with $t_2>t_1$, higher $L$ should increase the \SR. This is a consequence of the model; for if $t_2\gg t_1\gg t^*$, then $\frac{1}{t_2-t^*}\ll\frac{1}{t_1-t^*}$, and this may in turn be much smaller than $\mu$ or $k\mu$. This higher value of $\mu$ or $k\mu$ removing individuals by death would effectively swamp the effect of removal due to retirement. For given time $L$ that individuals are in the population, increasing $s_0\in[\frac{1}{3},\frac{2}{3}]$ increases $\mu$ and thus should tend the \SR\ towards $1$. Naturally, increasing $k$ should decrease the \SR\ by selectively removing males earlier.

\subsection{Numerical methods}
\label{sec:numerics}
We explore our model in two ways, both using Matlab's \verb|ode15s| function to integrate system~\eqref{eqn:DEsystem} numerically, ensuring that all variables remain non-negative. The first exploration is to straightforwardly study the parameter landscape on grids varying two life history parameters ($L$ and $t_1$), with other parameters fixed at various choices. The second is a broader sensitivity analysis varying all parameters (as well as initial condition) within reasonable ranges, allowing us to assess how each parameter influences \SR\ and dominant strategy.

\section{Results}
\label{sec:results}
We generate a grid on the parameter space and examine the resulting \SR\ and strategy contours in the plane whose axes are $L$ and $t_1$ (the $Lt_1$-plane), illustrated in Figures~\ref{fig:contours1}-\ref{fig:contours4}. We show several slices of the grid for different paternity theft and pair-bond break-up rates (assessing the behaviour with regard to effectiveness of guarding). The four figures show the results for all combinations of $t_2\in\{60,80\}$ and $k\in\{1,1.2\}$. Each figure shows a grid of six contour plots (landscapes) for paternity theft rate $g = 0$ or $0.003$, and pair-bond break-up rates $\beta = 0$, $\frac{1}{16}$, and $\frac{1}{4}$. Initial conditions for Figures~\ref{fig:contours1}-\ref{fig:contours4} had $G(0) = 250$ adult searching guarding males, $M(0) = 250$ adult multiple-mating males, $F(0) = 500$ receptive adult females, $F^G(0) = 0$ guarded females, $F^M(0) = 0$ unreceptive females, $C^G(0) = 500(1-g)$ guarded juveniles, and $C^M(0) = 500(1+g)$ unguarded juveniles (accounting for paternity theft in the juvenile population). Further simulations were done with variations in the initial male population structure, to assess whether or not initial conditions could influence the equilibrium. This is discussed in subsection~\ref{sec:bistability}.

Each landscape has an ``extinction boundary'' along the bottom and left of the image, where either menopause or death occurs at such a low average age that the average number of offspring per female is below replacement levels with the given mean birth rate $\rho = \frac{1}{3}$. In general, there is a clear boundary between regions of the $Lt_1$ plane where one strategy becomes completely prevalent and the other is suppressed; except near the extinction boundary to the lower left in some of the landscapes, this boundary has a positive gradient.

Contours of constant \SR\ are generally of positive gradient also. This observation is to be expected, as with increasing age female fertility ends $t_1$, the number of fertile females relative to fertile males increases; and as life expectancy at birth increases, the adult death rate $\mu$ decreases, causing females and male retirement rates $\tau$ and $\lambda$ to dominate in the removal of fertile adults from the population, increasing the \SR\ when $t_1<t_2$. Rephrasing the preceding logic, high adult mortality swamps removal by retirement; conversely, low adult mortality makes visible any differences between $t_1$ and $t_2$. In these plots, we do not have $t_1>t_2$, so only positive-gradient contours are visible; had we explored either lower $t_2$ or higher $t_1$, negative-gradient contours would appear.

There is a tipping point where either life expectancy at birth $L$ increases sufficiently or age female fertility ends $t_1$ decreases sufficiently that guarding takes over from multiple mating. Increasing either the pair-bond break-up rate $\beta$ or the paternity theft rate $g$ increases the proportion of the landscape that is taken up by multiple mating.

\subsection{Effects of male-specific mortality and retirement}
\label{sec:maleSpecificEffects}
\begin{figure}\centering
\begin{subfigure}[b]{4cm}
	\includegraphics[width=\textwidth]{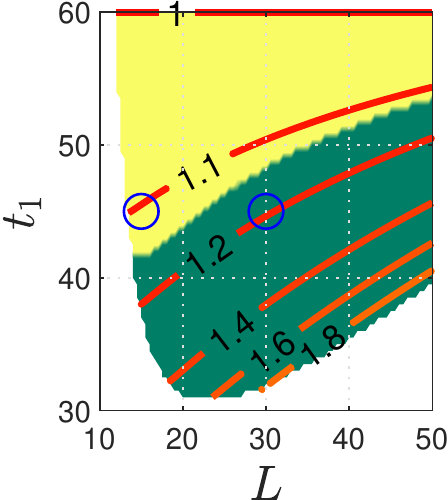}
	\caption{$g = 0$, $\beta = 0$}
	\label{subfig:contours1:1:1}
\end{subfigure}%
\begin{subfigure}[b]{4cm}
	\includegraphics[width=\textwidth]{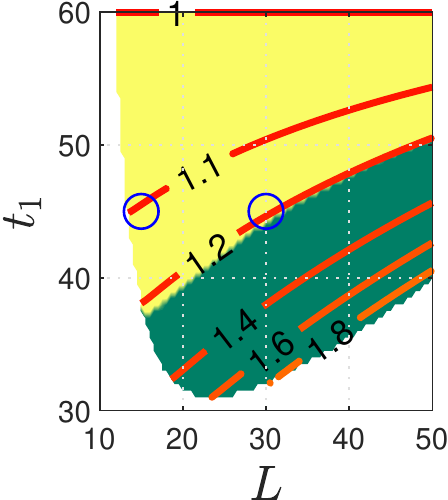}
	\caption{$g = 0$, $\beta = 1/16$}
	\label{subfig:contours1:2:1}
\end{subfigure}%
\begin{subfigure}[b]{4cm}
	\includegraphics[width=\textwidth]{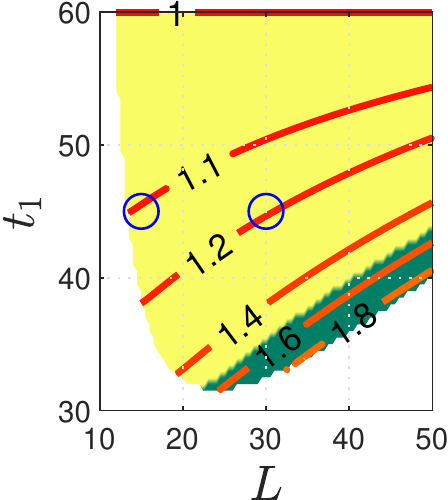}
	\caption{$g = 0$, $\beta = 1/4$}
	\label{subfig:contours1:3:1}
\end{subfigure}\\[1em]%
\begin{subfigure}[b]{4cm}
	\includegraphics[width=\textwidth]{fig_2_1_1_1_1}
	\caption{$g = 1/10$, $\beta = 0$}
	\label{subfig:contours1:1:2}
\end{subfigure}%
\begin{subfigure}[b]{4cm}
	\includegraphics[width=\textwidth]{fig_2_1_1_2_1}
	\caption{$g = 1/10$, $\beta = 1/16$}
	\label{subfig:contours1:2:2}
\end{subfigure}%
\begin{subfigure}[b]{4cm}
	\includegraphics[width=\textwidth]{fig_2_1_1_3_1}
	\caption{$g = 1/10$, $\beta = 1/4$}
	\label{subfig:contours1:3:2}
\end{subfigure}
\caption{Contour plots in the $Lt_1$-plane (longevity--age female fertility ends) of equilibrium values of \SR\ (numbered contours in red) and dominant strategy (yellow/green regions) as outcomes of numerically integrating equations~\eqref{eqn:DEsystem}. Initial fraction of males guarding $0.5$, male retirement age $t_2 = 60$, adult male relative mortality risk $k = 1.0$, with proportion surviving until maturity $s_0 = 1/2$, mean birth rate $\rho = 1/3$, carrying capacity $\nu = 1/1000$, couple-forming rate $r = 2$, and return rate for females mated with multiple maters $\sigma = 1$. Subplots compare outcomes for different values of pair-bond break-up rate $\beta = 0$, $1/16$, $1/4$ (columns) and paternity theft rate $g = 0$, $1/10$ (rows). White background indicates extinctions, dark/green indicates guarding dominates, light/yellow indicates multiple mating dominates. Labelled contours indicate contours of constants \SR. Small circles indicate approximate locations of chimpanzee and hunter-gatherer life history values at $(L,t_1)\approx(15,45)$, $(30,45)$ respectively.}
\label{fig:contours1}
\end{figure}

Figure~\ref{fig:contours1} shows landscapes for age of male retirement $t_2=60$, male-specific mortality risk $k=1$. Observe that because the male-specific mortality risk $k=1$ and age of male infirmity $t_2=60$, along the line with age female fertility ends $t_1 = 60$ we have \SR\ exactly equal to $1$, as we should expect. We see in this set of landscapes that the contour along which strategy switches tends to stay very close to contours of constant \SR.

\begin{figure}\centering
\begin{subfigure}[b]{4cm}
	\includegraphics[width=\textwidth]{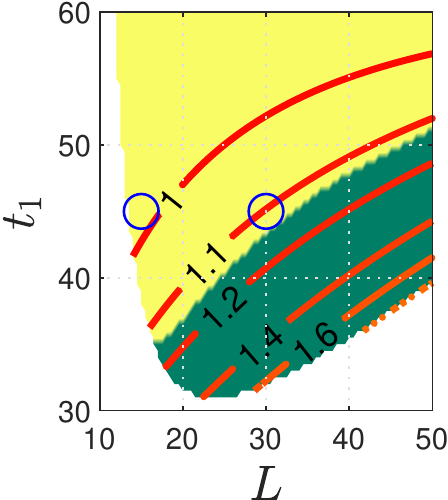}
	\caption{$g = 0$, $\beta = 0$}
	\label{subfig:contours2:1:1}
\end{subfigure}%
\begin{subfigure}[b]{4cm}
	\includegraphics[width=\textwidth]{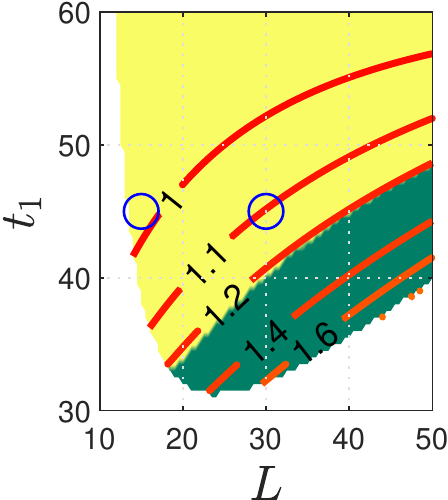}
	\caption{$g = 0$, $\beta = 1/16$}
	\label{subfig:contours2:2:1}
\end{subfigure}%
\begin{subfigure}[b]{4cm}
	\includegraphics[width=\textwidth]{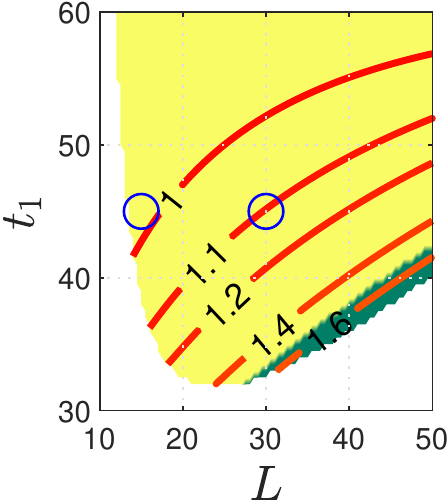}
	\caption{$g = 0$, $\beta = 1/4$}
	\label{subfig:contours2:3:1}
\end{subfigure}\\[1em]%
\begin{subfigure}[b]{4cm}
	\includegraphics[width=\textwidth]{fig_2_1_2_1_1}
	\caption{$g = 1/10$, $\beta = 0$}
	\label{subfig:contours2:1:2}
\end{subfigure}%
\begin{subfigure}[b]{4cm}
	\includegraphics[width=\textwidth]{fig_2_1_2_2_1}
	\caption{$g = 1/10$, $\beta = 1/16$}
	\label{subfig:contours2:2:2}
\end{subfigure}%
\begin{subfigure}[b]{4cm}
	\includegraphics[width=\textwidth]{fig_2_1_2_3_1}
	\caption{$g = 1/10$, $\beta = 1/4$}
	\label{subfig:contours2:3:2}
\end{subfigure}
\caption{Results of simulations as in Figure~\ref{fig:contours1}, now with male retirement age $t_2 = 60$, adult male relative mortality risk $k = 1.2$.}
\label{fig:contours2}
\end{figure}

In Figure~\ref{fig:contours2}, the age of male infirmity remains $60$, but male-specific mortality risk $k=1.2$. This is accompanied by an increase in the approximate \SR\ at which strategy appears to switch. While again the switch in strategy does not lie strictly \emph{on} a contour of constant \SR, it remains very near---arguably nearer than in Figure~\ref{fig:contours1}.

\begin{figure}\centering
\begin{subfigure}[b]{4cm}
	\includegraphics[width=\textwidth]{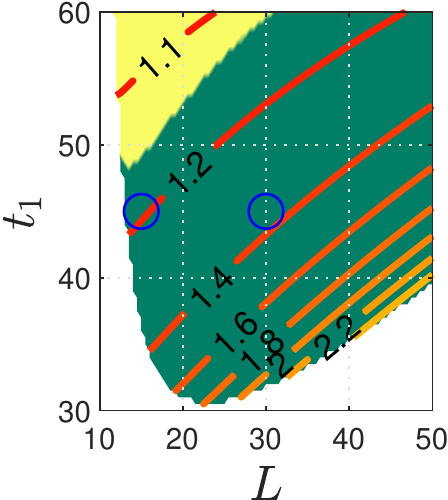}
	\caption{$g = 0$, $\beta = 0$}
	\label{subfig:contours3:1:1}
\end{subfigure}%
\begin{subfigure}[b]{4cm}
	\includegraphics[width=\textwidth]{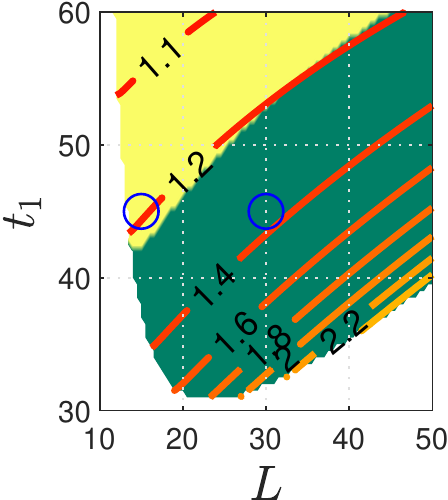}
	\caption{$g = 0$, $\beta = 1/16$}
	\label{subfig:contours3:2:1}
\end{subfigure}%
\begin{subfigure}[b]{4cm}
	\includegraphics[width=\textwidth]{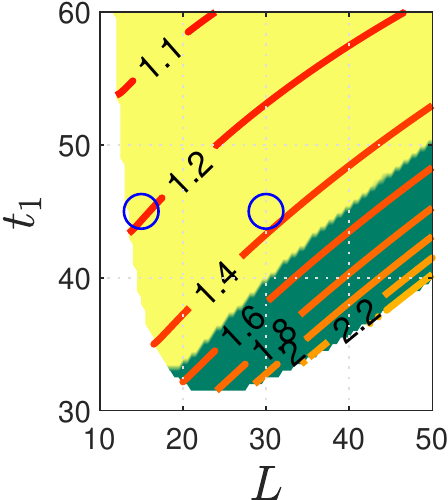}
	\caption{$g = 0$, $\beta = 1/4$}
	\label{subfig:contours3:3:1}
\end{subfigure}\\[1em]%
\begin{subfigure}[b]{4cm}
	\includegraphics[width=\textwidth]{fig_2_2_1_1_1}
	\caption{$g = 1/10$, $\beta = 0$}
	\label{subfig:contours3:1:2}
\end{subfigure}%
\begin{subfigure}[b]{4cm}
	\includegraphics[width=\textwidth]{fig_2_2_1_2_1}
	\caption{$g = 1/10$, $\beta = 1/16$}
	\label{subfig:contours3:2:2}
\end{subfigure}%
\begin{subfigure}[b]{4cm}
	\includegraphics[width=\textwidth]{fig_2_2_1_3_1}
	\caption{$g = 1/10$, $\beta = 1/4$}
	\label{subfig:contours3:3:2}
\end{subfigure}
\caption{Results of simulations as in Figure~\ref{fig:contours1}, now with male retirement age $t_2 = 80$, adult male relative mortality risk $k = 1.0$.}
\label{fig:contours3}
\end{figure}

Figure~\ref{fig:contours3} shows the same set of landscapes for $t_2 = 80$ and $k = 1$. The general features we observed earlier are largely the same, with the following exception: for $g = 0.003$, the contour along which strategy switches deviates markedly from the contours of constant \SR\ for low $L$ and low $t_1$, near the extinction boundary. If observed closely, this detail may be seen in the corresponding row of Figure~\ref{fig:contours1}. This, we posit, is an artefact of our life history model: for sufficiently low $L$, adult life expectancy may decrease more slowly with $L$, leading to comparatively longer adult life expectancies, with the upshot that males can therefore expect a relatively longer time to compete for the comparatively quickly retiring females, and to gain paternities when pair-bonded.

\begin{figure}\centering
\begin{subfigure}[b]{4cm}
	\includegraphics[width=\textwidth]{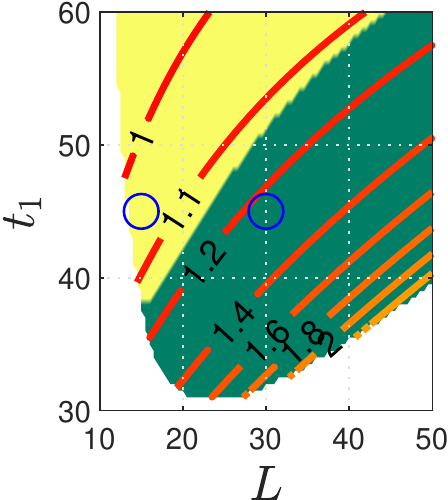}
	\caption{$g = 0$, $\beta = 0$}
	\label{subfig:contours4:1:1}
\end{subfigure}%
\begin{subfigure}[b]{4cm}
	\includegraphics[width=\textwidth]{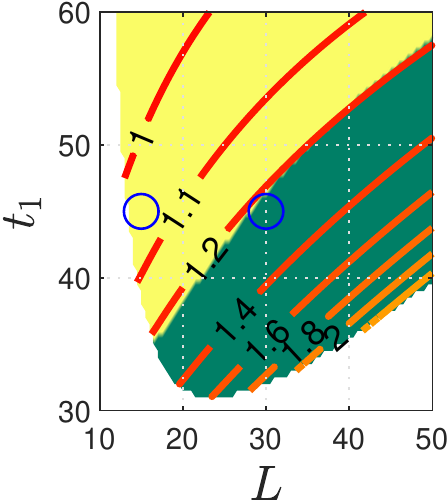}
	\caption{$g = 0$, $\beta = 1/16$}
	\label{subfig:contours4:2:1}
\end{subfigure}%
\begin{subfigure}[b]{4cm}
	\includegraphics[width=\textwidth]{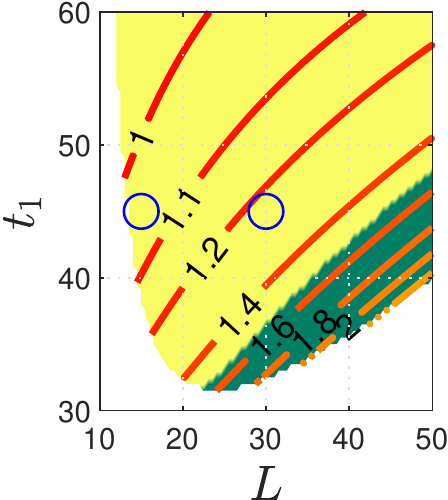}
	\caption{$g = 0$, $\beta = 1/4$}
	\label{subfig:contours4:3:1}
\end{subfigure}\\[1em]%
\begin{subfigure}[b]{4cm}
	\includegraphics[width=\textwidth]{fig_2_2_2_1_1}
	\caption{$g = 1/10$, $\beta = 0$}
	\label{subfig:contours4:1:2}
\end{subfigure}%
\begin{subfigure}[b]{4cm}
	\includegraphics[width=\textwidth]{fig_2_2_2_2_1}
	\caption{$g = 1/10$, $\beta = 1/16$}
	\label{subfig:contours4:2:2}
\end{subfigure}%
\begin{subfigure}[b]{4cm}
	\includegraphics[width=\textwidth]{fig_2_2_2_3_1}
	\caption{$g = 1/10$, $\beta = 1/4$}
	\label{subfig:contours4:3:2}
\end{subfigure}
\caption{Results of simulations as in Figure~\ref{fig:contours1}, now with male retirement age $t_2 = 80$, adult male relative mortality risk $k = 1.2$.}
\label{fig:contours4}
\end{figure}

Finally, Figure~\ref{fig:contours4} has $t_2 = 80$ and $k = 1.2$. Because the male-specific mortality risk factor is higher than in Figure~\ref{fig:contours3}, multiple mating gains more traction in this instance; however, because the age of male infirmity is higher than it is in Figure~\ref{fig:contours2}, guarding has more traction than in that instance. Again, note that the strategy-switch contour lies very close to contours of constant \SR.

The broad results indicate that in a landscape with varying life expectancy $L$ and age female fertility ends $t_1$ four things are of note:
\begin{enumerate}[1)]
\item multiple mating dominates when the average length of time that males can compete decreases (either by increasing the male-specific mortality risk $k$ or decreasing the age of male retirement $t_2$); \label{paralist:first}
\item \SR\ increases when male retirement is later in life, and when male-specific mortality is lower;
\item improving the effectiveness of guarding increases the region of the $Lt_1$-plane where guarding dominates; and
\item the boundary between multiple mating and guarding aligns well, overall, with contours of constant \SR\ in the $Lt_1$-plane, indicating that the ``signal'' of dominant male strategy as a function of \SR\ does not depend very much on life history parameters, but does depend on behavioural parameters such as paternity theft rate and the average length of time a pair bond lasts.
\end{enumerate}
Point~\ref{paralist:first} above indicates that given a fixed effectiveness of guarding, the length of time that males are able to compete for paternities determines how they are best able to maximise their average number of paternities. Less time means spreading risk across multiple mates, and more time means investing in a single mate.

Assuming that all other parameters are constant, it may be possible to shift from multiple mating to guarding as the dominant male strategy just by increasing life expectancy at birth, a manoeuvre that simultaneously increases the \SR\ when $t_1<t_2$. A plausible instance of this idea would be to draw a path connecting the two circles in the upper-middle plot of Figure~\ref{fig:contours4}, which has pair-bond break-up rate $\beta=\frac{1}{16}$ and paternity theft rate $g=0$. One might imagine that our common ancestor with chimpanzees (as close as can be represented in a model such as this) lived near the point with $L=15$ and $t_1=45$, and a ``short'' evolutionary trajectory that increased only our life span (as far as this model's parameters are concerned, i.e. with no increase in age that female fertility ends) could have tipped us from multiple mating into guarding male mating strategies.

\subsection{A region of bistability}
\label{sec:bistability}
\begin{figure}\centering
\begin{subfigure}[b]{4cm}
	\includegraphics[width=\textwidth]{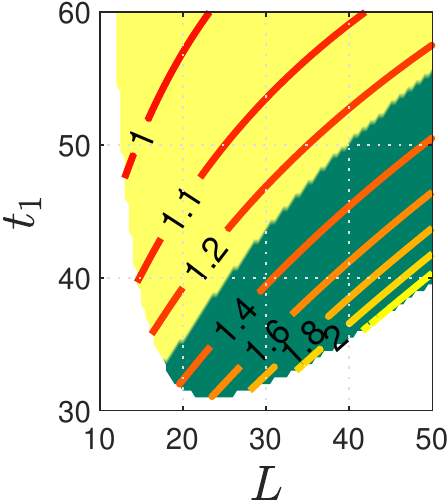}
	\caption{$R_0 = 0.1$}
	\label{subfig:bistability1}
\end{subfigure}%
\begin{subfigure}[b]{4cm}
	\includegraphics[width=\textwidth]{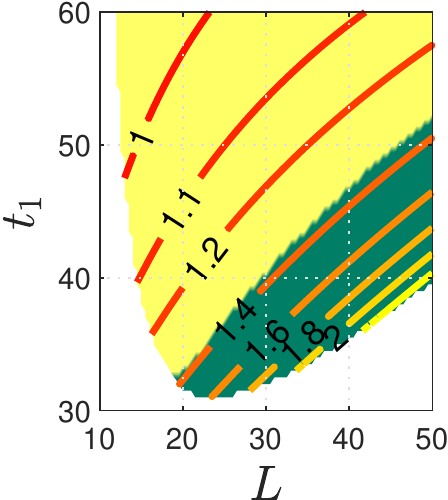}
	\caption{$R_0 = 0.9$}
	\label{subfig:bistability2}
\end{subfigure}
\caption{Contour plots similar to subplots of Figures~\ref{fig:contours1}-\ref{fig:contours4}, showing the effect of different initial population structures when paternity theft is in effect ($g = 0.1$). Left~\ref{subfig:bistability1}, initial fraction of adult population multiple-mating $R_0 = 0.1$; right~\ref{subfig:bistability2}, initial fraction of adult male multiple-mating $R_0 = 0.9$. Observe that when more multiple-mating males are initially present (right) the yellow area corresponding to multiple-mating dominating is larger, close to the contour with \SR\ $\approx1.4$. Parameters are otherwise as in lower-middle subplot of Figure~\ref{fig:contours4}, with pair-bond break-up rate $\beta = 1/16$, male retirement age $t_2 = 80$, adult male relative mortality risk $k = 1.2$, with proportion surviving until maturity $s_0 = 1/2$, mean birth rate $\rho = 1/3$, carrying capacity $\nu = 1/1000$, couple-forming rate $r = 2$, and return rate for females mated with multiple maters $\sigma = 1$.}
\label{fig:initialPopulationEffect}
\end{figure}

Suppose an initial population contains a mixture of strategies---that is, the populations $G(0)+F^G(0)>0$ and $M(0)>0$, and the initial fraction of multiple-maters is $R_0 = M(0)/(M(0)+G(0)+F^G(0))$. If $g=0$ and no paternity theft occurs we naturally expect the equilibrium populations to not depend on $R_0$ (which is indeed so). However, if $g>0$ is this still the case? Simulations with different values of $R_0$ and paternity theft rate $g>0$ show that indeed the initial population structure affects the population structure at equilibrium. We illustrate this in Figure~\ref{fig:initialPopulationEffect}, where the initial populations have $R_0 = 0.1$(left) and $R_0 = 0.9$ (right), with $g = 0.1$. Juvenile populations are in proportion to the adult male structure, with an adjustment of $(1\pm g)$ for guarded ($-$) and unguarded ($-$) juvenile populations.

The contour along which the dominant strategy changes when $R_0 = 0.1$ lies approximately between \SR\ values of about $1.2$ and $1.3$, whereas when $R_0 = 0.9$ this contour is very close to the contour with \SR\ approximately $1.4$. These observations indicate a region of bistability, for which a future bifurcation analysis could prove interesting; however, for low levels of paternity theft the region of bistability appears to be small, even if not non-negligible.

\subsection{Sensitivity analysis}
\label{subsec:sensitivityAnalysis}
We use Latin hypercube sampling to generate $100\,000$ parameter points within reasonable ranges. We take age female fertility ends $t_1$ between $40$ and $55$, age of male retirement $t_2$ between $60$ and $80$, birth frequency $\rho$ between $\frac{1}{4cm}$ and $\frac{1}{2.5}$, and male relative risk $k$ between $1$ and $1.2$, and excluding cases where the resulting \SR\ was outside the range of $\frac{1}{3}$ to $3$. All other parameter ranges are as specified in Table~\ref{tab:params}. Results are given shown in Table~\ref{tab:spearmanAll} for the $78\,226$ points that satisfy the \SR\ requirement. The Spearman partial rank correlation coefficients are $\rho_\text{\SR}$ for the change in \SR\ due to each variable, and similarly $\rho_R$ for the change in the fraction $R:= M/(M+G+F^G)$ of multiple maters in the fertile adult male population. We classify a correlation as ``very weak'', ``weak'', ``moderate'', ``strong'', or ``very strong'' if the absolute value of the correlation coefficient is in the respective ranges $[0,0.2)$, $[0.2,0.4)$, $[0.4,0.6)$, $[0.6,0.8)$, or $[0.8,1]$. Dominant strategy is determined by $R$ at equilibrium: if $R<0.5$ guarding is dominant, and if $R>0.5$ multiple-mating is dominant.

Of particular note are the following observations:
\begin{enumerate}
\item \label{enum:SRresults} \SR\ is most strongly positively correlated with age of male retirement $t_2$ and mean longevity $L$ (in descending order), as may be expected from the parameter ranges, for which $t_2>t_1$, so increasing either of these two parameters will result in relatively more males. Similarly, \SR\ is most strongly negatively correlated (in decreasing order of correlation strength) with age of female retirement $t_1$, birth rate $\rho$, male-specific mortality risk $k$, and proportion $s_0$ surviving to mean age of maturity. All other parameters have a very weak correlation with \SR, $|\rho_\text{\SR}|<0.1$.
\item \label{enum:stratResults1} Age of female retirement $t_1$ has the strongest correlation with dominant strategy; in particular, increasing age of female retirement (for the parameter range we've chosen, this means bringing it closer to parity with the age of male retirement) correlates with an increase in the dominance of multiple-mating. Other life history parameters have much lower correlation with dominant strategy.
\item \label{enum:stratResults2} Outside the life history parameters, the return rate $\sigma$ of unreceptive females $F^M$ has the strongest (but only a moderate) correlation with dominant strategy, similar, but opposite in sign, to pair-bond break-up rate $\beta$.
\item \label{enum:stratResults3} Paternity theft rate $g$ is moderately correlated with strategy for the range tested, as is birth rate $\rho$. All other parameters have very weak correlations with dominant strategy.
\end{enumerate}

\begin{table}
\centering
\caption{Spearman partial rank correlation coefficients $\rho_\text{\SR}$ and $\rho_R$ for, respectively, \SR\ and fraction $R = M/(M+G+F^G)$ of males multiple-mating. Due to the number of samples satisfying $1/3<\text{\SR}<3$ ($\sim8\times10^4$), most $p$-values were $<0.05$. Exceptions (marked with~$^*$) were paternity theft rate $g$ and $\rho_{\text{\SR}}$ ($p=0.34201$), pair-bond break-up rate $\beta$ and $\rho_{\text{\SR}}$ ($p=0.18465$), and initial proportion of males multiple-mating $R_0 = M(0)/(M(0)+G(0)+F^G(0))$ and $\rho_{\text{\SR}}$ ($p=0.73475$). The implication is that the correlation in these cases could not be distinguished from no average rank difference due to changes in each of those respective variables.}
\label{tab:spearmanAll}
\begin{tabular}{rSSSS}
Variable	& $\rho_\text{\SR}$	& $\rho_R$	\\
\hline
$L      $	&  0.31598			& -0.096241			\\
$s_0    $	& -0.13266			&  0.10749			\\
$\delta $	&  0.017005			& -0.022914			\\
$\mu    $	& -0.083075			&  0.025083			\\
$t_1    $	& -0.95766			&  0.72761			\\
$t_2    $	&  0.67168			& -0.23224			\\
$\rho   $	& -0.42046			&  0.46645			\\
$\nu    $	&  0.010335			& -0.17244			\\
$r      $	& -0.0081298		&  0.039327			\\
$g		$	& -0.0033977$^*$	&  0.5508			\\
$\beta	$	&  0.0047435$^*$	&  0.46109			\\
$\sigma	$	&  0.0095936		& -0.5626			\\
$k		$	& -0.57252			&  0.18527			\\
$R_0	$	& -0.0012115$^*$	&  0.056362			\\
\hline
\end{tabular}
\end{table}

\subsection{Discussion}
\label{subsec:discussion}
Our sensitivity analysis, in Subsection~\ref{subsec:sensitivityAnalysis}, shows that the \SR\ is controlled largely by the parameters one would expect: ages $t_1$ and $t_2$ of female and male retirement, and male-specific mortality risk $k$, as well as, to some extent, the mean longevity $L$. Variation in $L$ can in a certain sense ``expose'' the difference between $t_1$ and $t_2$; if $L$ is very low, we would expect individuals to die before having a chance to be removed from the fertile pool due to infirmity or menopause. The birth rate $\rho$ may affect \SR\ by causing a population to more quickly reach its carrying capacity (determined by $\nu$). Since males are removed by death according to $\mu(k+\nu P)$, this may result in an extra discrepancy than if they were removed by a different death law (such as $k\mu(1+\nu P)$). Allowing $k<1$ may have also consequently resulted in a less monotonic relationship between birth rate and \SR.

The correlation of birth rate $\rho$ to strategy is unsurprising: a higher birth rate implies more opportunities for paternity \emph{and} paternity theft \citep{hawkesRobertsCharnov1995}, hence increasing $\rho$ tends to increase the frequency of multiple mating as dominant strategy. The likely underlying cause is the fact that multiple maters can recruit multiple females at once; increasing the average number of offspring per recruitment is therefore a significant advantage if you do not restrict yourself to a single female as guarders do.

In chimpanzees, the inter-birth interval is approximately four to five years, but in humans it is on average less than three---very surprising for a species of our longevity. The Grandmother Hypothesis posits that grandmothering subsidies for their daughters' child-rearing shortened it, and retaining the ancestral age of that female fertility ends. The difference in inter-birth interval thus suggests in our model that humans should possibly be more likely to multiple-mate than chimps, but the weakness of the correlation and differences in other parameters tending towards the success of guarding may be sufficient to compensate. A further difference is that human longevities are greater, while most male chimpanzees are dead by $t_1$, thus having a much shorter typical period in which to compete for paternities.

Similarly surprising is the comparative strength of the dependence of dominant strategy on $\nu$. Although very weakly correlated, the correlation is far from negligible. As $\frac{1}{\nu}$ corresponds to the carrying capacity, increasing the carrying capacity (decreasing $\nu$) appears to increase the success of multiple mating. We posit that this is because although the effect of changing $\nu$ on the \SR\ is minimal, the dominant strategy is strongly dependent on the number of receptive females at equilibrium. Consequently, it may be that as $\nu$ increases (corresponding to a decrease in the maximum population size), the number of receptive females decreases, meaning that fewer receptive females can be recruited by multiple maters, giving guarders an advantage as they obtain more assured paternities with what mates they can obtain and guard for an extended length of time.

Results indicated in Table~\ref{tab:spearmanAll} show that female and male retirement ages $t_1$ and $t_2$, birth rate $\rho$, and male-specific mortality risk $k$ have non-negligible correlations with both \SR\ and dominant strategy. Moreover, each of these parameters has an alternate relationship with each of those two quantities: a parameter that correlates with an increase in \SR\ correlates with a decrease in multiple-mating as dominant, or vice-versa. Most other parameters correlate to only \SR\ or only dominant strategy, else at most weakly to both. The correlation between dominant strategy and \SR\ has Spearman partial rank correlation $-0.6166$, controlling for $\delta$ and $\mu$ (child and adult death rates), which are dependent on the choice of longevity $L$ and proportion $s_0$ surviving at the mean age of maturity. Controlling in addition for $t_1$, $t_2$, $\rho$ and $k$, the correlation reduces to $-0.2824$. That is, there is only a weakly monotonic correlation between \SR\ and dominant strategy under variation of all other parameters. This result is consistent with the observations in Figures~\ref{fig:contours1}-\ref{fig:contours4}, where the contour along which dominant strategy switches is largely (but not exactly) aligned with some contour of constant \SR.

What we can understand is that there is no distinct value of \SR\ at which strategy switches in the full parameter space. \SR\ does, however, act as a guide: guarding does not appear to occur frequently for more balanced \SR s (approximately $1$). What we suggest guides the dynamic is that guarding becomes favourable when the female retirement time scale, $\frac{1}{\tau}$, is short enough that it becomes significant in causing retirements for females in the unreceptive $F^M$ pool. A female recruited by a guarding male will typically produce $\frac{\rho}{\beta+\tau+\mu}$ offspring, compared to $\frac{\rho}{\sigma+\tau+\mu}$ for one recruited by a multiple mater. Increasing $\tau$ (the rate at which females advance to menopause) not only reduces the average length of time that females spend in either reproducing pool, but will reduce the number of receptive females for all males to recruit. However, guarding males who manage to recruit a receptive female are guaranteed a longer time with them to obtain paternities, but the success of multiple mating depends on both how many receptive females females and also on how many competitors there are: fewer possible recruitments drastically reduces the possible number of paternities, and thus guarding is more likely to out-compete multiple mating.

\section{Conclusion}

We have developed a model for a primate-like population endowed with sufficient details about life history and behaviour to begin exploring mathematically the link between the \SR\ and the dominant strategy employed by males, and the elements of behaviour and biology that may affect the outcome to varying extents. We wish to use this model (or subsequent models) to help explain behavioural differences between ourselves and our nearest living relatives, the chimpanzees, with the benefit that the model is general enough that the approach employed herein can be used to understand the mating strategy dynamics of other species as well. Our model does not assume anything about the details of male-male or male-female interactions---only their net outcomes as they relate to mean reproduction and mortality rates, though it is a detail specific to humans to distinguish the age that female fertility ends $t_1<t_2$ the age that male fertility ends.

Our results indicate that female scarcity relative to the length of the male's window of opportunity to reproduce strongly dictates the likelihood that guarding will take hold over multiple mating, given some understanding of how mate-guarding is accomplished. That is, given an estimate of the effectiveness of guarding (both in maintaining long term guarded relationships and in preventing paternity theft) and, in relevant species (such as humans) the lengths of the reproductive intervals of males and females, the \SR\ may provide an index by which to predict the typical mating strategy of males. Particularly, a greater abundance of fertile males relative to potentially fertile females indicates that a guarding strategy is more likely to be employed by males.

Though we obtain interesting and possibly useful results through this model, there remain some limitations and some questions arise. For instance, direct patrilineal inheritance of strategy is arguably unrealistic; it would be interesting to understand to what extent a propensity to use a particular strategy is inherited versus flexible and adaptable to different social and cultural situations within an individual male's lifetime. Parallel to that: to what extent does guarding behaviour really exclude seeking paternities outside the pair bond? Additionally, how does female choice affect outcomes, and what form does it take? Dawkins' famous ``Battle of the sexes'' game \citep{dawkins1976}, re-examined, for instance, in \citep{mcnamaraEtAl2009}, frames female choice in terms of being ``coy'' versus ``fast'', but the question of strategy depends also on other aspects of male behaviour: as mentioned in the introduction, male chimpanzees can be violently coercive with females and dangerous to their offspring, so the latter must make choices in light of that knowledge to ensure their own safety and the safety of their young. Finally, what other strategies might males employ that could affect the balance that our model predicts? It is possible that guarding males also benefit their offspring through provision of care and protection from infanticide, phenomena we chose not to model here, and even that the promise of such care (as in the Battle of the sexes) makes a male more attractive as a mate and thus more likely to obtain (or retain) a mate.

\section*{Acknowledgments}
DR and PSK were supported by the Australian Research Council, Discovery Project (DP160101597). The authors also wish to extend their thanks to the editor and reviewers for their thorough evaluation of our work, significantly improving its final form.

\bibliography{bibliography}

\end{document}